\documentclass[twoside]{irmaems}
\usepackage{amssymb} 
\usepackage{amsmath} 
\usepackage{latexsym}

\renewcommand\theequation{\arabic{section}.\arabic{equation}}

\theoremstyle{definition} 

\theoremstyle{plain}      

\input epsf

\begin{document}

\title{ 3-Dimensional Schlaefli Formula and Its
Generalization  }

\author{Feng Luo\thanks{ Work partially supported by the NSF.} }

\address{
Department of Mathematics\\
Rutgers University\\
Piscataway, NJ 08854, USA\\
email:\,\tt{ fluo\@math.rutgers.edu}\\
\\
\\
\text{Dedicated to the memory of Xiao-Song Lin} } \maketitle

\begin{abstract}
Several identities similar to the Schlaefli formula are
established for tetrahedra in a space of constant curvature.

\end{abstract}

\tableofcontents

\section{Introduction}\label{s-1}


 One of the most important identities in
low-dimensional geometry is the Schlaefli formula. It states that
for a tetrahedron in a constant curvature $\lambda= \pm 1$ space,
the volume $V$, the length $x_{ij}$, and the dihedral angle
$a_{ij}$ at the ij-th edge are related by

\begin{equation}\label{1.1}
\frac{ \partial V}{\partial a_{ij}} = \frac{\lambda}{2} x_{ij}
\end{equation}
where $V = V(a_{12}, a_{13}, a_{14}, a_{23}, a_{24}, a_{34})$ is a
function of the angles. See for instance \cite{Mil} or \cite{AV}
for a proof.

\medskip

\epsfxsize=2truein\centerline{\epsfbox{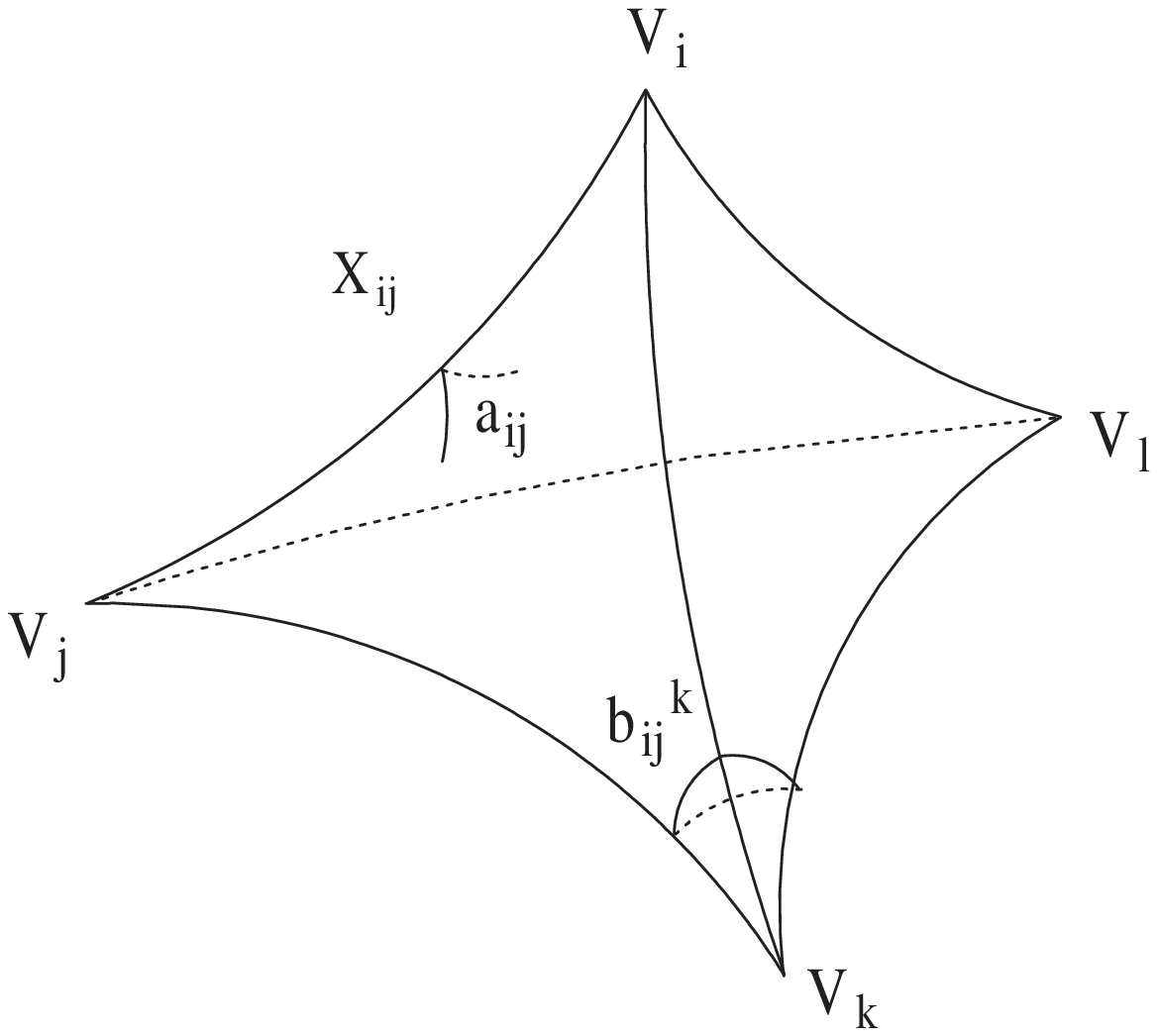}}

\centerline{Figure 1.1}

In 3-dimensional polyhedral geometry, a space is obtained by
isometrically gluing tetrahedra along their codimension-1 faces.
The metric is determined by the edge lengths and the curvature at
an edge is $2\pi$ less the sum of dihedral angles at the edge.
From this point of view, the Schlaefli formula relates the most
important geometric quantities: the volume, the metric (=edge
lengths) and the curvature (= dihedral angles) in a simple elegant
identity. The Schlaefli formula plays a vital role in a
variational principle for triangulated 3-manifolds. See for
instance Regge's work \cite{Re} on discrete general relativity.

One consequence of (1.1) is that differential 1-forms
\begin{equation}\label{1.2}
\sum x_{ij} da_{ij} \quad  \text{and} \quad \sum a_{ij} dx_{ij}
\end{equation} are closed.

We may recover the volume function  $V$ in (1.1) by integrating
the 1-form $\sum x_{ij} da_{ij}$. Thus the closeness of the
1-forms in (1.2) essentially captures the Schlaefli formula.

The basic problem in polyhedral geometry is to understand the
relationship between the metric and its curvature. In the case of
tetrahedra, this prompts us to study the curvature map $K(x)=a$
sending the edge length $x=(x_{12}, x_{13}, x_{14}, x_{23},
x_{24}, x_{34})$ to the dihedral angle $a=(a_{12}, a_{13}, a_{14},
a_{23}, a_{24}, a_{34})$.  The Jacobian matrix $D(K)$ of the
curvature map is  the $6 \times 6$ matrix $[\frac{\partial
a_{ij}}{\partial x_{rs}}]_{6 \times 6}$. The closeness of the
1-forms in (1.2) is equivalent to say that the Jacobian matrix
$D(K)$ is symmetric.  It turns out the Jacobian matrix
$[\frac{\partial a_{ij}}{\partial x_{rs}}]_{6 \times 6}$ enjoys
many more symmetries. One of the symmetry was discovered by E.
Wigner \cite{Wi} and Taylor-Woodward \cite{TW}. The purpose of
this paper is to find the complete set of symmetries of the
Jacobian matrix $D(K)$ of the curvature map. These symmetries
should have applications in 3-dimensional topology and geometry.
In particular, the relationships between the Jacobian matrix
$D(K)$, the 6j symbols, the quantum 6j symbols, and the volume
conjecture are very attractive problems. See for instance the work
of \cite{TW} and \cite{Rob}.

The complete set of symmetries was discovered by us a few years
ago. We thank Walter Neumann for suggesting us to write it up for
publication.

This paper is dedicated to the memory of Xiao-Song Lin who made
important contributions to low-dimensional topology. He was a
great colleague and friend.

The paper is organized as follows. In \S2, we state the main
theorem. These theorems are proved in \S3. A more general version
of it involving complex valued lengths and angles can be found in
\cite{Lu1}.

\section{The main theorem}

 Let a tetrahedron in
$\bf S^3$ or $\bf H^3$ or $\bf E^3$ have vertices $v_1,
v_2,v_3,v_4$. Let $a_{ij}$ and $x_{ij}$ be the dihedral angle and
the edge length at the ij-th edge $v_i v_j$. We consider the angle
$a_{ij}$ as a function of the lengths $x_{12}, x_{13}, x_{14},
x_{23}, x_{24}, x_{34}$.

Here is our main result.

\bigskip
\noindent {\bf Theorem 2.1.} \it Define  $P_{rs}^{ij} =
\frac{1}{\sin(a_{ij}) \sin(a_{rs})}\frac{\partial a_{ij}}{\partial
x_{rs}} $. Then these quantities satisfy the following identities
for any tetrahedron in any of the three geometries $\bf S^3, \bf
H^3, \bf E^3$. The indices $i,j,k,l$ are assumed to be pairwise
distinct.

(i) (Schlaefli) $P^{ij}_{rs} = P_{ij}^{rs}$.

(ii) (\cite{Wi}, \cite{TW}) $P^{ij}_{kl} =P^{ik}_{jl} =
P^{il}_{jk}$.

(iii) $P^{ij}_{ik} = -P^{ij}_{kl} \cos a_{jk}$.

(iv) $P^{ij}_{ij}= P^{ij}_{kl} w_{ij}$ where
 $w_{ij} =(c_{ij}c_{jk}c_{ki} + c_{ij}c_{jl}c_{li} +
c_{ik}c_{jl} + c_{il}c_{jk})/\sin^2(a_{ij})$ and $c_{rs} =
\cos(a_{rs})$.

(v) $P^{ij}_{rs} = P^{i'j'}_{r's'}$ where $\{i,j\} \neq \{r,s\}$
and for a subset $\{a, b\} \subset \{i,j,k,l\}$, the set $\{a',
b'\}$ is $\{i,j,k,l\}-\{a,b\}$. \rm


\bigskip
In the spaces $\bf S^3$ and $\bf H^3$ of constant curvature
$\lambda =\pm 1$, a tetrahedron is determined by its dihedral
angles $a_{ij}$. Thus the length $x_{ij}$ can be considered as a
function of the angles. The similar theorem is,

\bigskip \noindent
 {\bf Theorem 2.2.} \it Define  $R_{rs}^{ij} = \frac{1}{\sin(\sqrt{\lambda}x_{ij}) \sin(\sqrt{\lambda}x_{rs})}
\frac{\partial x_{ij}}{\partial a_{rs} }$. Then these quantities
satisfy the following identities for any tetrahedron in spherical
and hyperbolic geometries. Let the indices $i,j,k,l$ be distinct.

(i) (Schlaefli) $R^{ij}_{rs} = R_{ij}^{rs}$.

(ii) (\cite{Wi}, \cite{TW}) $R^{ij}_{kl} =R^{ik}_{jl} =
R^{il}_{jk}$.

(iii) $R^{ij}_{ik} = R^{ij}_{kl} \cos (\sqrt{\lambda} x_{il})$.

(iv) $R^{ij}_{ij}= R^{ij}_{kl} w_{ij}$ where
$$w_{ij}= \frac{- c_{ij}c_{ik}c_{il} -c_{ji}c_{jk}c_{jl} + c_{ik}c_{jl} + c_{il} c_{jk}}
{ \sin^2( \sqrt{\lambda} x_{ij})},$$ and
 $c_{rs} =\cos(\sqrt{\lambda} x_{rs})$,

(v) $R^{ij}_{rs} = R^{i'j'}_{r's'}$ where $\{i,j\} \neq \{r,s\}$
and for a subset $\{a, b\} \subset \{i,j,k,l\}$, the set $\{a',
b'\}$ is $\{i,j,k,l\}-\{a,b\}$.

\rm

\medskip
We remark that the matrices $[\frac{\partial a_{ij}}{\partial
x_{rs}}]$ and $[\frac{\partial x_{ij}}{\partial a_{rs}}]$ are
inverse of each other when $\lambda \neq 0$. Theorem 2.2  follows
from theorem 2.1 by taking the dual.  Indeed, in the spherical
tetrahedral case, the dual tetrahedron has dihedral angle
$\pi-x_{ij}$ and edge length $\pi-a_{ij}$ at the kl-th edge of the
dual simplex. Thus, theorem 2.2 follows. The hyperbolic tetrahedra
case in theorem 2.2 can be deduced from spherical case by
analytical continuation.

Theorem 2.1 suggests that the matrix $M D(K) M=[P^{ij}_{rs}]_{6
\times 6}$ where $M$ is the diagonal matrix whose diagonal entries
are $\frac{1}{\sin( a_{ij})}$ exhibits more symmetries than the
Jacobian matrix $D(K)$.

Both theorems are special cases of a complex valued edge-length
and dihedral angle relation. This will be discussed in \cite{Lu1}.

\medskip
\section{A proof of theorem 2.1}

\medskip

We need to recall the cosine law and its derivative form in order
to compute the Jacobian matrix $[\frac{\partial a_{ij}}{\partial
x_{rs}}]$ effectively.

Let $K^2=$ $\bf S^2$, or $\bf H^2$ or $\bf E^2$ be the space of
constant curvature $\lambda =1,-1$, or $ 0$. Define a function
$S_{\lambda}(t)$ as follows. $S_{0}(t) =t$; $S_1(t) = \sin(t)$ and
$S_{-1}(t) =\sinh(t)$. The sine law for a triangle of lengths
$l_1, l_2, l_3$ and opposite angles $a_1, a_2, a_3$ in $K^2$ can
be stated as
\begin{equation}\label{Sine Law}
\frac{S_{\lambda}(l_i)}{\sin(a_i)}=\frac{S_{\lambda}(l_j)}{\sin(a_j)}
\end{equation}
A different way to state the sine law is that the expression
$$ A_{ijk} = \sin(a_i) S_{\lambda}(l_j)S_{\lambda}(l_k)$$
is symmetric in indices $i,j,k$ where $\{i,j,k\}=\{1,2,3\}$. For
this reason, we call $A_{ijk}=A_{123}$ the \it A-invariant \rm of
the triangle.

\medskip
\noindent {\bf Proposition 3.1.}(\cite{CL}, \cite{Lu}) \it  Let a
triangle in $K^2$ have inner angles $a_1, a_2, a_3$ and edge
lengths $l_1, l_2, l_3$ so that $l_i$-th edge is opposite to the
angle $a_i$. Then

(i) $\frac{\partial a_i}{\partial l_j} =-\frac{\partial
a_i}{\partial l_i} \cos( a_k)$ where $\{i,j,k\}=\{1,2,3\}$,

(ii) $\frac{\partial a_i}{\partial l_i} =
\frac{S_{\lambda}(l_i)}{A_{123}}$ \rm

\medskip
See \cite{CL} or \cite{Lu} for a proof.
\bigskip

Let us introduce some notations before beginning the proof. The
indices $i,j,k,l$ are pairwise distinct, i.e.,
$\{i,j,k,l\}=\{1,2,3,4\}$. The face triangle $\Delta v_i v_j v_k$
will be denoted $\Delta ijk$. The inner angle at the vertex $v_k$
of the triangle $\Delta ijk$ is denoted by $b^k_{ij}$. The link at
the vertex $v_k$, denoted by $Lk(v_k)$, is a spherical triangle
with edge lengths $b^k_{ji}, b^k_{il}, b^k_{lj}$ and inner angles
$a_{ki}, a_{kj}, a_{kl}$ so that $a_{ki}$ is opposite to
$b^k_{jl}$. The A-invariant of the triangle $\Delta ijk$ is
denoted by $A_{ijk}$.

In the calculation below, we consider $b^i_{jk}$ as a function of
$x_{rs}$'s using the cosine law for the triangle $\Delta ijk$. By
the definition, we have, \begin{equation} \label{1} \frac{\partial
b^i_{jk}}{\partial x_{rs}}=0 \end{equation} if $\{r,s\}$ is not a
subset of $ \{i,j,k\}$.
 The function $a_{ij}$ is considered as a
function of $b^r_{st}$'s by the cosine law applied to either the
link Lk($v_i$) or Lk($v_j$). In this way the dihedral angle
$a_{ij}$, when considered as a function of the lengths $x_{rs}$'s,
is a composition function.

To prove theorem 2.1, note that identity (i) in theorem 2.1 is the
Schlaefli formula (1.2). Identity (v) follows from identity (iii).
By symmetry, we only need to consider three partial derivatives:
$\frac{\partial a_{ij}}{\partial x_{kl}}$, $\frac{\partial
a_{ij}}{\partial x_{ik}}$ and $\frac{\partial a_{ij}}{\partial
x_{ij}}$.


\bigskip
\subsection{The partial derivatives $\frac{\partial a_{ij}}{\partial x_{ik}}$ and $\frac{\partial a_{ij}}{\partial
x_{kl}}$}

Consider the link Lk($v_i$). Using proposition 3.1(ii), the chain
rule and  (3.2), we have (see Fig. 2.1(a)),
\begin{equation}\label{2} \frac{\partial a_{ij}}{\partial x_{kl}} = \frac{\partial a_{ij}}{\partial
b^i_{kl}}\frac{  \partial b^i_{kl}}{\partial x_{kl}} =
\frac{\partial a_{ij}}{\partial b^i_{kl}} \frac{
S_{\lambda}(x_{kl})}{A_{ikl}}.
\end{equation}

Similarly, using Lk($v_j$), we have
\begin{equation}\label{2} \frac{\partial a_{ij}}{\partial x_{kl}} = \frac{\partial a_{ij}}{\partial
b^j_{kl}}\frac{  \partial b^j_{kl}}{\partial x_{kl}} =
\frac{\partial a_{ij}}{\partial b^j_{kl}} \frac{
S_{\lambda}(x_{kl})}{A_{jkl}}.
\end{equation}

\medskip

\epsfxsize=3truein\centerline{\epsfbox{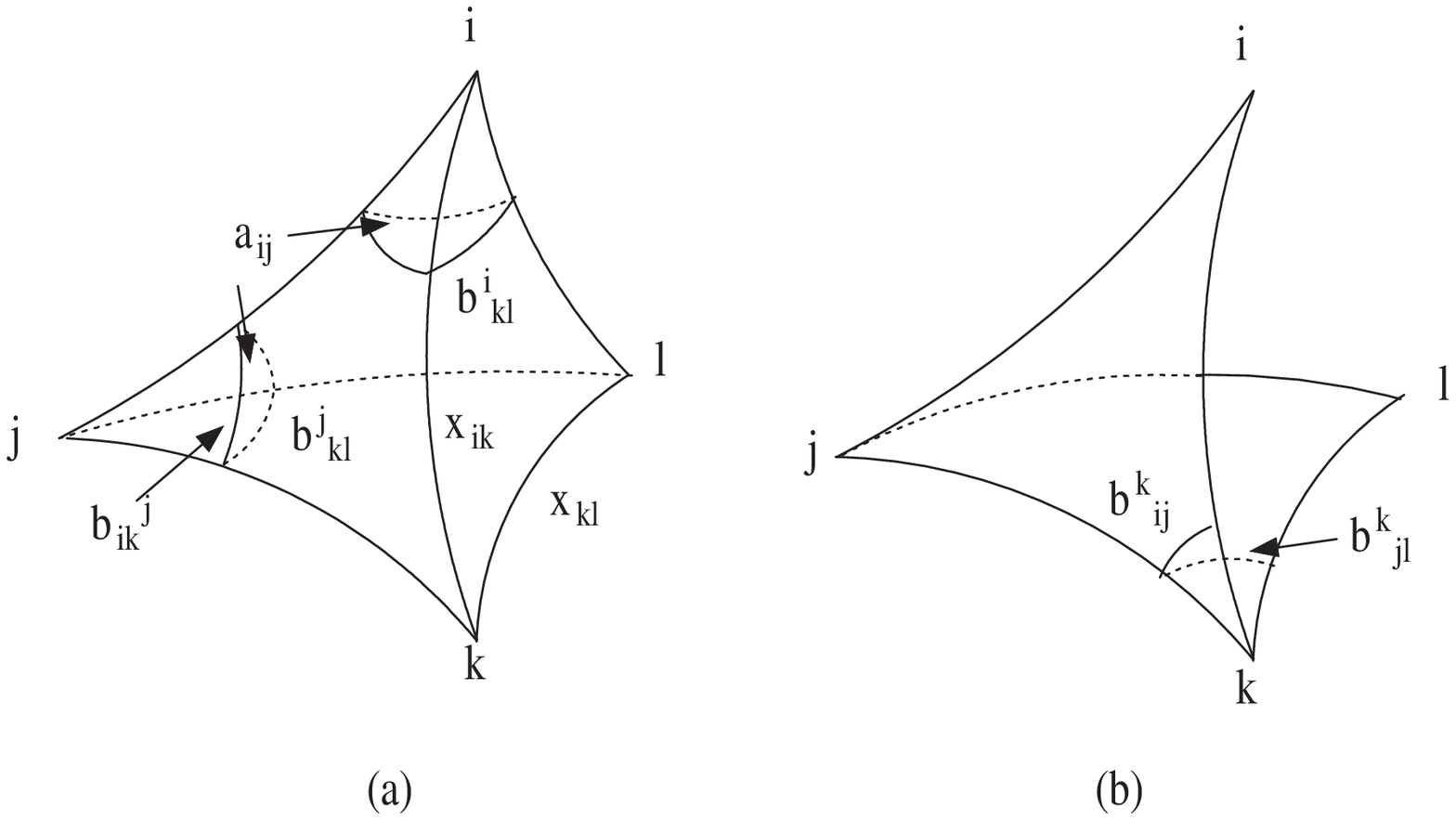}}

\centerline{Figure 2.1}

\medskip
Now we use the link Lk($v_j$) to find $\frac{\partial
a_{ij}}{\partial x_{ik}}$. By (3.2) and the chain rule, we have
\begin{equation}\label{a} \frac{\partial a_{ij}}{\partial x_{ik}}
= \frac{\partial a_{ij}}{\partial b^j_{ik}} \frac{
\partial b^j_{ik}}{\partial x_{ik}}. \end{equation} By
proposition 3.1 applied to Lk($v_j)$ and $\Delta ijk$, we see
(3.5) is equal to,
\begin{equation} \label{a}
-\frac{ \partial a_{ij}}{\partial b^j_{kl}} \cos(a_{jk})
\frac{S_{\lambda}(x_{ik})}{A_{ijk}}.
\end{equation}
Using (3.4), we can write (3.6) as,
\begin{equation} \label{}
 -\frac{\partial a_{ij}}{\partial x_{kl}} \cos(a_{jk})
\frac{A_{jkl}S_{\lambda}(x_{ik})}{A_{ijk} S_{\lambda}(x_{kl})}.
\end{equation}

Now by the definition of the A-invariant of triangles $\Delta ijk$
and $\Delta jkl$ (see Fig. 2.1(b)), we have,
\begin{equation} A_{jkl} =
S_{\lambda}(x_{jk})S_{\lambda}(x_{kl})\sin(b^k_{jl}) \quad
\text{and} \quad A_{ijk} =
S_{\lambda}(x_{jk})S_{\lambda}(x_{ik})\sin(b^k_{ij}).\end{equation}
Thus (3.7) can be simplified to
\begin{equation}\label{a} -\frac{\partial a_{ij}}{\partial x_{kl}} \cos(a_{jk})\frac{
\sin(b^k_{jl})}{\sin(b^k_{ij})}.
\end{equation}
By the sine law applied to the spherical triangle Lk($v_k$), we
see (3.9) is equal to
\begin{equation} \label {q} \frac{\partial a_{ij}}{\partial x_{ik}} = -\frac{\partial a_{ij}}{\partial
x_{kl}} \cos(a_{jk}) \frac{\sin(a_{ik})}{\sin(a_{kl})}.
\end{equation}

This is equivalent to identity (iii),
\begin{equation} \label{q} P^{ij}_{ik} = -P^{ij}_{kl} \cos a_{jk}.
\end{equation}

Use the Schlaefli formula that $P^{ij}_{ik} =P^{ik}_{ij}$, we
obtain from (3.11)

$$-P^{ij}_{kl} \cos(a_{jk}) = -P^{ik}_{jl} \cos(a_{jk}).$$

This shows that \begin{equation}\label{1} P^{ij}_{kl} =
P^{ik}_{jl}.\end{equation} By symmetry, identity (ii) holds for
all indices.


\bigskip
\subsection{The partial derivative  $\frac{\partial a_{ij}}{\partial
x_{ij}}$}

By (3.2), the chain rule, we have, in the triangle Lk($v_i$),

\begin{equation}\label{1} \frac{\partial a_{ij}}{\partial x_{ij}} =
\frac{\partial a_{ij}}{\partial b^i_{jk}} \frac{ \partial
b^i_{jk}}{\partial x_{ij}} + \frac{
\partial a_{ij}}{\partial b^i_{jl}} \frac{ \partial b^i_{jl}}{\partial x_{ij}}. \end{equation}

 Using proposition 3.1, we see that (3.13) is equal to
\begin{equation} \label{ 2} \frac{ \partial a_{ij}}{\partial b^i_{kl}} \cos(a_{ik}) \cos(b^j_{ik})\frac{ S_{\lambda}(x_{jk})}{A_{ijk}} +
\frac{ \partial a_{ij}}{\partial b^i_{kl}} \cos(a_{il})
\cos(b^j_{il}) \frac{S_{\lambda}(x_{jl})}{A_{ijl}}. \end{equation}

Using (3.3), we see (3.14) is equal to \begin{equation} \label{2}
 \frac{\partial a_{ij}}{\partial x_{kl}} [ \cos (a_{ik}) \cos
 (b^j_{ik})
\frac{S_{\lambda}(x_{jk})A_{ikl}}{S_{\lambda}(x_{kl})A_{ijk}} +
\cos(a_{il})\cos(b^j_{il})
\frac{S_{\lambda}(x_{jl})A_{ikl}}{S_{\lambda}(x_{kl}) A_{ijl}}].
\end{equation}

Using the sine law for triangles $\Delta ikl$ and $\Delta ijl$ as
in (3.8), we can rewrite (3.15) as
\begin{equation}\label{2}
\frac{\partial a_{ij}}{\partial x_{kl}} [ \cos(a_{ik})
\cos(b^j_{ik}) \frac{ \sin(b^k_{il})}{\sin(b^k_{ij})} +
 \cos(a_{il}) \cos(b^j_{il}) \frac{\sin(b^l_{ik})}{\sin(b^l_{ij})}]. \end{equation}

 Using the sine law in triangles Lk($v_k$) and Lk($v_l$), we see that (3.16)
 is the same as
 \begin{equation}
 \frac{\partial a_{ij}}{\partial x_{kl}} [ \cos(a_{ik}) \cos(b^j_{ik}) \frac{\sin(a_{kj})}{\sin(a_{kl})} +
\cos(a_{il})\cos(b^j_{il})\frac{ \sin(a_{lj})}{\sin(a_{lk})}].
\end{equation}

\begin{equation}
=P^{ij}_{kl} [ \cos(a_{ik}) \cos(b^j_{ik})
\sin(a_{kj})\sin(a_{ij}) + \cos(a_{il})\cos(b^j_{il})
\sin(a_{lj})\sin(a_{ij})].
\end{equation}

On the other hand, by the cosine  law for the spherical triangle
Lk($v_j$), we have
$$
\cos(b^j_{ik}) \sin(a_{kj})\sin(a_{ij}) = \cos a_{kj} \cos a_{ij}
+ \cos a_{lj}. $$ and

$$\cos(b^j_{il})
\sin(a_{lj})\sin(a_{ij}) = \cos a_{jl} \cos a_{ij} + \cos
a_{jk}.$$

Substitute these into (3.18),  we obtain

$$\frac{ \partial a_{ij}}{\partial x_{ij}} =P^{ij}_{kl} [c_{ij}c_{jk}c_{ki} +c_{ij}c_{jl}c_{li}
+c_{ik}c_{jl} + c_{il}c_{jk}]
$$ where $c_{rs} =\cos(a_{rs})$.

This is the identity (iv) since $P^{ij}_{ij} =
\frac{1}{\sin^2(a_{ij})} \frac{\partial a_{ij}}{\partial x_{ij}}$.

\end{document}